\documentclass[12pt,a4paper,oneside]{article}
\usepackage{amsmath,amstext,amssymb,amsopn,amsthm,eufrak}
\usepackage[latin1]{inputenc}

\theoremstyle{plain}
\newtheorem{thm}{Theorem}[section]

\newtheorem{lem}[thm]{Lemma}

\theoremstyle{definition}

\theoremstyle{remark}
\newtheorem*{rem*}{Remark}
\newtheorem{rem}{Remark}

\newcommand{\R}{\mathbb{R}}
\newcommand{\N}{\mathbb{N}}
\newcommand{\Q}{\mathbb{Q}}

\renewcommand{\leq}{\leqslant}
\renewcommand{\geq}{\geqslant}

\title { On regularity properties of Bessel flow }
\date{}
\author{Vostrikova L.\\ LAREMA, D\'epartement de Math\'ematiques,
Universit\'e d'Angers, France}

\begin{document}
\maketitle

\begin{abstract}
We study the differentiability of Bessel flow $\rho : x \rightarrow
\rho ^x_t$, where $(\rho ^x_t)_{t\geq 0}$ is BES $^x(\delta $) process of
dimension $\delta >1$ starting from $x$. For $\delta \geq 2$ we prove the existence of
bicontinuous derivatives in P-a.s. sense at $x\geq 0$ and we study the
asymptotic behaviour of the derivatives at $x=0$. For $1< \delta <2$
we prove the existence of a modification of Bessel flow having
derivatives in probability sense at $x\geq 0$. We study the asymptotic
behaviour of the derivatives at $t=\tau_0(x)$ where $\tau _0(x)$ is
the first zero of $(\rho ^x_t)_{t\geq 0}$.  
\end{abstract}
\footnotetext{
\emph{2000 MS Classification:} 60GXX, 60G17  
   \\
  {\it Key words and phrases}:  Bessel flow, regularity, bicontinuity
  . \\}

 %%%%%%%%%%%%%%%%%%%%%%%%%%%%%%%%%%%%%%%%%%%%%%%%%%%%%%%%%%%%%%%%%%%%%%%%%%%%%%
                         \section{Introduction}
 %%%%%%%%%%%%%%%%%%%%%%%%%%%%%%%%%%%%%%%%%%%%%%%%%%%%%%%%%%%%%%%%%%%%%%%%%%%%%%
 
 The regularity of flows of diffusion processes is an important problem 
 related to the stability of solutions of SDEs with respect 
 to the initial value. This problem is well-studied when the 
 coefficients of the diffusion equation are regular( cf.{\sc Kunita }\cite{KU},
 {\sc Protter}\cite{P}). Some results for the non-Lipschitz case is
 given in {\sc Ren, Zhang }\cite{RZ}.
 
 As is well-known,  Bessel squared process of dimension $\delta > 0$, denoted by 
 $\mbox{BESQ}^{x^2}(\delta )$, starting from $x^2$, is the unique strong solution of the following 
 stochastic differential equation: for all $x\geq 0$ and $t>0$
 \begin{equation}\label{besq}
 X_t^x = x^2 + 2\displaystyle\int_0^t\sqrt{X_s^x} d\beta _s + \delta t,
 \end{equation}
 where $\beta = (\beta _t)_{t\geq 0}$ is standard Brownian motion. 
 For $x>\gamma >0$ with fixed $\gamma $,  $X^x = (X_t^x)_{t\geq 0}$ is 
 diffusion process with locally Lipschitz coefficients on 
 $(0, +\infty )$. Moreover, the 
 derivatives of diffusion coefficients with respect to initial value  
 are also locally Lipschitz on the same set. It gives, using the comparison theorem,  
 that the flow of $\mbox{BESQ}^{x^2}(\delta )$ processes with $ x > \gamma ,$ is a diffeomorphisme
 up to explosion time for derivatives, which is 
 $$ \tau _0(\gamma ) = \inf\{t\geq 0 : X^{\gamma }_t = 0\}$$
 where $\inf\{\emptyset \} = \infty $.
 We remark that in the case $\delta \geq 2$ we have  $P(\tau _0(\gamma ) = \infty) = 1$ , and 
 in the case $1<\delta <2$ we get $P(\tau _0(\gamma ) < \infty) = 1$. In general, we 
 cannot expect to establish some regularity properties  after explosion time for the derivatives.
 But $\mbox{BESQ}^{x^2}(\delta )$ process is a very special case in which this study
 may be possible. It should be noticed that being particular, BESQ
 processes appear relatively often: it is so for radial part squared
 of Brownian motion; the laws of some local times for  Brownian
 motion are related to BESQ process, and the same is true for some
 processes related with running maximum of Brownian
 motion (see {\sc Borodin, Salminen} \cite{BS}). Another important
 example is given by the trace of Wishart process which is also
 BESQ process (see {\sc Bru}\cite{B}). In this context  we should also mention Dunkl
 process which  radial part squared is also BESQ process ( see {\sc
 L. Gallardo, M. Yor}\cite{GY1},\cite {GY2}).

 The same comments can be made for $\mbox{BES}^x(\delta )$  with $\delta >0$. This
 process is related to $\mbox{BESQ}^{x^2}(\delta )$ in the following
 way: for $t\geq 0$
 $$\rho ^x_t = \sqrt{X_t^{x}}.$$
 It is well known (see, for instance \cite{RY}, chapter XI) that for $\delta >1$, the $\mbox{BES}^x(\delta )$
 process is the solution of the following differential  equation: for all $x\geq 0$ and $t>0$
 \begin{equation}\label{bes}
 \rho ^x_t = x + \beta _t + \displaystyle\frac{(\delta -1)}{2}
 \displaystyle\int_0^t\displaystyle\frac{1}{\rho ^x_s}ds
 \end{equation}
 where $\beta = (\beta _t)_{t\geq 0}$ is standard Brownian motion.
 For $\delta =1$ the $\mbox{BES}^x(\delta )$ process satisfies: for all $x\geq
 0$ and $t>0$
 \begin{equation}\label{bes1}
 \rho ^x_t = x + \beta _t + L^0_t(x)
 \end{equation}
 where $L^0_t(x)$ is local time of $\mbox{BES}^x(\delta )$ process at
 zero.
 For $0<\delta <1$ $\mbox{BES}^x(\delta )$ process verify: for all $x\geq
 0$ and $t>0$
 \begin{equation}\label{bes2}
 \rho ^x_t = x + \beta _t +  \displaystyle\frac{(\delta -1)}{2}
 \mbox{v.p.}\displaystyle\int_0^t\displaystyle\frac{1}{\rho ^x_s}ds
 \end{equation}
 where the integral in the right-hand side is understanding  is in
 v.p. sense.
  
 The structure of $\mbox{BES}^x(\delta )$ process is simpler then the
 one of $\mbox{BESQ}^{x^2}(\delta )$ 
 process in a sense that the equations (\ref{bes}), (\ref{bes1}), (\ref{bes2})
 do not contain a stochastic integral. This is
 the reason why we focuss our study on the flow of $\mbox{BES}^x(\delta )$ processes. We remark 
 that some indications related to the regularity property of Bessel flow
 with $\delta \geq 2$ and $x>0$ can be found in 
 {\sc Hirsch, Song}\cite {HS}.
   
 The aim of this paper is to study the regularity property of the flow of
 $\mbox{BES}^x(\delta )$ processes  with $\delta >1$.  We will distinguish
 the cases of $\delta \geq 2$ and $1<\delta <2$, and inside of them also 
 the cases $x>0$  and $x\geq 0$.
 
 \vskip 0.5cm

 %%%%%%%%%%%%%%%%%%%
\begin{thm}\label{t1}
For $\delta \geq 2$, the flow of $\mbox{BES}^x(\delta )$ processes 
has (P-a.s.) derivatives of all orders  with respect to $x$ for $x>0$ which are 
bicontinuous in $(x,t)$ on the set $]0, +\infty [\times [0, +\infty[$.
\end{thm}
\vskip 0.5cm 
To prove Theorem \ref{t1}  we reduce first the problem to the case of $x>\gamma >0$,
then we do localisation and we use the classical results.
The case $x\geq 0$ is very different from the  case $x>0$ from point
of view of properties and, then also from technical point of view.  
In the case $x\geq 0$ the mentionned above procedure does not work and we have to use some
identity in law and some fine asymptotics of Spitzer type to conclude
(cf. {\sc Spitzer}\cite{SP}, {\sc Messulam P., Yor M.}\cite{MY}).

\begin{thm}\label{t11}
 For $\delta = 2$, the flow of $\mbox{BES}^x(\delta )$ processes has
 (P-a.s.) derivatives of all orders n  at $x=~0$ (and a fortiori for
 $x>0$). These derivatives are bicontinuous in $(x,t)$ on the set 
$[0, +\infty [\times ]0, +\infty [.$
 Moreover, for the derivatives of Bessel flow $(\rho^x _t)_{t\geq 0, x\geq 0}$
 the following asymptotic relations hold.
\begin{enumerate}
\item[a)] For $n\geq 1$ uniformly on  compacts in $t$ of $(0, + \infty )$
 and P-a.s.
$$\displaystyle\lim _ {x\rightarrow 0+}\left(\ln\, (\,\displaystyle\frac{\partial ^n \rho ^x_t}{\partial x
  ^n}\,)\,/ \ln x \right) = +\infty.$$
\item[b)] For  $n\geq 1$ the convergence in law sense holds:
$$\displaystyle\lim_ {x\rightarrow 0+}\left(\ln\, (\,\displaystyle\frac{\partial ^n \rho ^x_t}{\partial x
  ^n}\,)\,/(\ln x)^2\right) = - T_1(\beta)/2$$
where $T_1(\beta )$ is the first passage time of the level 1 for
standard Brownian motion,
\item[c)]  Uniformly on  compacts in $t$ of
$(0, + \infty )$ and P-a.s.
$$\displaystyle\lim_ {x\rightarrow 0+}\left( x^{n-1}\,\,\displaystyle\frac{\partial ^n \rho ^x_t}{\partial x
  ^n}/\frac{\partial \rho ^x_t}{\partial x }\right)= U_{n-1}$$
where $U_{n-1} = U_1(U_1-1)\cdots (U_1-n+2)$ and $U_1$ is a random
variable given by (\ref{c1}) with $\nu = 1$.
\item[d)] For  $T> \epsilon >0$ and $0< \gamma < 1/(n-1)$ with the
  same  $U_{n-1}$ as in c)  
$$\displaystyle\lim_ {x\rightarrow 0+}E\left(\sup_{\epsilon \leq t\leq
  T}x^{n-1}\,\,
|\displaystyle\frac{\partial ^n \rho ^x_t}{\partial x
  ^n}/\frac{\partial \rho ^x_t}{\partial x }|\right)^{\gamma}=E(|
U_{n-1}|^{\gamma })$$
\end{enumerate}

For  $\delta > 2$, the flow of $\mbox{BES}^x(\delta )$ processes has (P-a.s.) at $x=0$ 
the derivatives only up to the order $n < n(\delta )$ where $n(\delta )=2 +
\displaystyle\frac{1}{\delta -2}$ ( and of all orders for $x>0$).
These derivatives are bicontinuous in $(x,t)$ on the set 
$[0, +\infty [\times ]0, +\infty [.$ Moreover  for $n< n(\delta )$  
we have:
\begin{enumerate}
\item[ a')] uniformly on  compacts in $t$ of $(0, + \infty )$
 and P-a.s.
$$\displaystyle\lim _ {x\rightarrow 0+}\left( \ln\, (\,\displaystyle\frac{\partial ^n \rho ^x_t}{\partial x
  ^n}\,)\,/ \ln x \right) = n(\delta ) -n,$$
\end{enumerate}
 and also the property c) and the property d) with $0< \gamma < \nu /(n-1)$ and
 $\nu =2\delta -3$ in (\ref{c1}).
\end{thm}
%%%%%%%%%%%%%%%%%%%%%
\vskip 0.5cm

\begin{rem}
For the regularity at $x=0$ we have the following picture. If $\delta \geq 3$ then
the flow has only two derivatives in
P-a.s. sense and no derivatives of order $n>2$ even in probability sense.
If $m\in \N ^*$
and $2+ \displaystyle\frac{1}{m+1}\leq \delta <2+ \displaystyle\frac{1}{m}$, then
the flow of $\mbox{BES}^x(\delta )$ processes has exactly $2+m$ derivatives in
P-a.s. sense. We remark that the regularity of the flow is increasing as 
$\delta \downarrow 2$,  and for $\delta = 2$ the flow is $C^{\infty }$.
The asymptotic relations a), b) give us logarithmic asymptotics for
n-th derivative of $\rho ^x_t$. The asymptotic relations c),
d) characterize the behaviour of the ratio of the n-th and the first
derivatives in P-a.s. and $L^{\gamma}$ - sense.

\end{rem} 
\vskip 0.5cm
 
If $1 < \delta < 2$, then $\mbox{BES}^x(\delta )$ process 
touches 0 with probability 1 and the results will be 
different from the previous case. To present the results
let us denote as before for $x>0$
\begin{equation}
\label{tauo}
\tau _0(x) = \inf \{ s \geq 0 : \rho^x _s = 0\}
\end{equation}
with $\inf\{\emptyset \} = +\infty $.
\vskip 0.5cm

%%%%%%%%%%%%%%%%%%%%%%%%%%
\begin{thm}\label{t2}
In the case $1 < \delta < 2$ and $x>0$ there exists a modification
$\tilde{\rho }$ of the Bessel flow 
in the space $D(\R ^{+,*}, C(\R ^{+,*})$ with the following properties:
\begin{enumerate}

\item[a)] $\tilde{\rho}$ is bicontinuous P-a.s. and has bicontinuous derivatives
of all orders on the set\\
$]0, +\infty [\times [\![ 0, \tau _0(x)[\![$. These derivatives coincide with the ones of 
$(\rho^x_{t\wedge \tau _0(x)})_{x>0, t>0}$.

\item[b)] For each $(x,t)$ with $x>0, t>0$  $\tilde{\rho}$ has
 derivatives in probability sense only up to the order 
$n<n(\delta )$ with $n(\delta ) =  \displaystyle\frac{1}{2-\delta }$,
which are bicontinuous in probability. Moreover, for $n<n(\delta )$ we have:
$$P \displaystyle\lim _ {y\rightarrow x+}\left(
\ln \,(\,\displaystyle\frac{\partial ^n \rho ^y_{\tau_0(x)}}{\partial y
  ^n}\,)\,/ \ln \,(\,\rho ^y_{\tau _0(x)}\,)\,\right) = n(\delta )-n,$$
\item[c)] For $n< n(\delta)$  we have
$$P\displaystyle\lim_ {y\rightarrow x+}\left((\,\rho ^y_{\tau_0(x)}\,)^{n-1}\,\,\,\,
\displaystyle\frac{\partial ^n \rho ^y_{\tau _0(x)}}{\partial y ^n}/
\displaystyle\frac{\partial \rho ^y_{\tau _0(x)}}{\partial y}\right)
= \frac{U_{n-1}}{x^{n-1}}$$
where $U_{n-1}$ is the
same as in Theorem \ref{t11} and  $U_1$ is a random
variable given by (\ref{c1}) with $\nu = 5-2\delta$.
\item[d)] For $n< n(\delta)$ and  $0< \gamma < (5-2\delta )/(n-1)$ with the
same  $U_{n-1}$ and $U_1$ as in c)
$$\displaystyle\lim_ {y\rightarrow x+}E\left((\,\rho ^y_{\tau _0(x)}\,)^{n-1}\,\,\,\,
|\displaystyle\frac{\partial ^n \rho ^y_{\tau _0(x)}}{\partial y ^n}/
\displaystyle\frac{\partial \rho ^y_{\tau _0(x)}}{\partial y}|\right)^{\gamma}
=  \,\frac{E(\,| U_{n-1}|^{\gamma}\,)}{x^{\,\gamma (n-1)}}$$
\end{enumerate}
 In the case $1<\delta <2$ and $x\geq 0,\, t>0$ the mentionned above
 modification of Bessel flow has the same regularity as for $x>0,\ t>0$.

\end{thm}
\vskip 0.5cm
%%%%%%%%%%%%%%%%%%%%%%%%% 
 
\begin{rem}
For the regularity in probability sense for $x>0$ and $t>0$ we have the
following picture.
For $1<\delta \leq 3/2$ the flow has only one derivative in
probability sense. 
For $m\in \N^*$ and $2 -\displaystyle\frac{1}{m+1}<\delta
 \leq 2 - \displaystyle\frac{1}{m+2}$ the considered modification
has exactly $m+1$ derivatives in probability sense. We remark that the regularity in probability sense
is increasing to infinity as $\delta \uparrow 2$. The interpetation of
the asymptotic relations is the same as in Remark 1.
\end{rem} 
\vskip 0.5cm
  
%%%%%%%%%%%%%%%%%%%%%%%%%%%%%%%%%%%%%%%%%%%%%%%%%%%%%%%%%%%%%%%%%%%%%%%%%%%%%%%%%%%%%%%%%%
             \section{Regularity of Bessel flow for $\delta \geq 2$}
%%%%%%%%%%%%%%%%%%%%%%%%%%%%%%%%%%%%%%%%%%%%%%%%%%%%%%%%%%%%%%%%%%%%%%%%%%%%%%%%%%%%%%%%%%
 \vskip 0.5cm
 
 In the case $x>0$ we begin with some rather general Lemmas.
  
 \begin{lem}\label{l1}
 Let $\rho ^x$ be the strong unique (P-a.s.) solution of the equation (\ref{bes})
 with initial value $x$. If for each $\gamma > 0$ the flow of  $\rho ^x$
 with $x>\gamma $ is bicontinuous and has (P-a.s.) bicontinuous in $(x,t)$
 derivatives of all orders with
 respect to $x$ on the set $]\gamma, +\infty[\times [0, +\infty[$
 , then there exists an extension of the flow on the set  $]0, +\infty[\times [0, +\infty[$ having the
 same properties.
 \end{lem}
 \vskip 0.5cm
 
 \noindent {\it Proof\,\,}
A simple patching  with respect to $\gamma $ proves the result.$\Box $
\vskip 0.5cm

Let now $\gamma >0$ be fixed and $x>\gamma >0$. To localise the coefficients of the
equation (\ref{bes}) we take a bicontinuous version of $\rho ^x$
(see \cite{RY}, p.362 ). For $0< \epsilon < \gamma $
we put
\begin{equation}\label{te}
\tau _{\epsilon} = \inf\{ t\geq 0 : \rho ^{\gamma  }_t \leq \epsilon\},
\end{equation}
and
\begin{equation}\label{to}
\tau  = \inf\{ t\geq 0 : \rho ^{\gamma }_t = 0\},
\end{equation}
with $\inf \{ \emptyset \} = \infty $.
To simplify the notations and since $\gamma $ is fixed, we do not write that 
$\tau _{\epsilon }$ and $\tau $ depend on $\gamma $.
\vskip 0.5cm

\begin{lem}\label{l2}
Suppose for each $\epsilon >0$ there exists a bicontinuous version of the
flow of the process $(\rho ^x_{t\wedge \tau_{\epsilon}})_{t\geq 0}$ having (P-a.s.) 
bicontinuous in $(x,t)$ derivatives with respect to $x$ on the set
$]\gamma , +\infty [\times [\![ 0, \tau _{\epsilon}[\![$ then there exists
an extension of the flow  having the same
properties on the set $]\gamma , +\infty [\times [\![ 0, \tau [\![$.
\end{lem}
\vskip 0.5cm

\noindent {\it Proof\,\,}A simple patching with respect to $\epsilon $
 gives the result.$\Box $
\vskip 0.5cm

\noindent {\it Proof of Theorem} \ref{t1}
Using the Lemmas \ref{l1}, \ref{l2} and the fact that $P(\tau = \infty) = 1$, we reduce our study
to the process $(\rho^x_{t\wedge \tau_{\epsilon }})_{t\geq 0}$ 
with $x\in ]\gamma , +\infty [$ and $\gamma >0$, where $\tau _{\epsilon }$
is defined by (\ref{te}).
By comparison theorem (P-a.s.) for all $t\geq 0$ and $x>\gamma $
$$\rho^{x}_{t\wedge \tau _{\epsilon}} \geq \rho^{\gamma }_{t\wedge \tau _{\epsilon}}
\geq \epsilon .$$
 
For  $x,y \in ]\gamma , +\infty [$ we denote
\begin{equation}
\label{incr1}
Z^{y,x}_{t\wedge \tau_{\epsilon}} =
\displaystyle
\frac{\rho ^y_{t\wedge \tau_{\epsilon}}
 -\rho ^x_{t\wedge \tau_{\epsilon}}}{y -x}.
\end{equation}
From (\ref{bes}) we obtain the following linear equation:
\begin{equation}
\label{incr2}
Z^{y,x}_{t\wedge \tau_{\epsilon}} =
1 - \displaystyle\frac{(\delta -1)}{2}
\displaystyle\int_0^{t\wedge \tau_{\epsilon}}
\displaystyle\frac{Z^{y,x}_s}{\rho ^y_s\rho ^x_s}ds
\end{equation}
and, hence, the solution
\begin{equation}
\label{incr3}
Z^{y,x}_{t\wedge \tau_{\epsilon}} = 
\exp \left\{ - \displaystyle\frac{(\delta -1)}{2}
\displaystyle\int_0^{t\wedge \tau_{\epsilon}}
\displaystyle\frac{1}{\rho ^y_s\rho ^x_s}ds \right\}.
\end{equation}
To take the limit as $y\rightarrow x$ we use bicontinuity of the flow of 
$\rho ^x$, the fact that
on the interval $[\![0, t\wedge \tau_{\epsilon}[\![$ we have a minoration:
$\rho ^y_s\geq \epsilon ,\:\rho ^x_s\geq \epsilon $.
By Lebesgue dominating convergence theorem, the first 
derivative of the flow is given by:
\begin{equation}\label{d1}
Y^x_{t\wedge \tau_{\epsilon}} =
\exp \left\{ - \displaystyle\frac{(\delta -1)}{2}
\displaystyle\int_0^{t\wedge \tau_{\epsilon}}
\displaystyle\frac{1}{(\rho ^x_s)^2}ds \right\}.
\end{equation} 
 The bicontinuity of the first
derivative follows in the same way using bicontinuity of the flow $\rho ^x$
and the above minoration.

Now, using the expression of the first derivative and the arguments
mentionned for bicontinuity of the first derivative,  we prove, in recurrence way,
the existence and bicontinuity of the $n$-th derivative.
$\Box $
\vskip 0.5cm

%%%%%%%%%%%%%%%%%%%%%%%%%%%%%%case x=0%%%%%%%%%%%%%%%%%%%%%%%%%%%%

To study the existence and bicontinuity of the derivatives of Bessel
flow at $x=0$, we need a renforced scaling property and some asymptotic results.

\vskip 0.5cm

%%%%%%%%%%%%%%%%%%%%%
\begin{lem}\label{sl}
If we consider bicontinuous versions of Bessel processes, then for all $c>0$
the renforced scaling property holds:
$$\mathcal L 
\left( \displaystyle\frac{1}{c}(\rho ^x_{c^2t})_{t\geq 0, x>0}\right) =
\mathcal L 
\left( (\rho ^{x/c}_{t})_{t\geq 0, x>0}\right) .
$$
\end{lem}
\vskip 0.5cm

\noindent {\it Proof\:} This is clair result for a finite number of $x$, say
$x_1,x_2,\cdots ,x_n$, due to simple scaling property of Bessel
process and the uniqueness of the solution of (\ref{bes}). Then, the
result follows by continuity. $\Box $ 

%%%%%%%%%%%%%%%%%%%%%
\begin{rem}
\label{r1a}
As a corollary of this Lemma, the law of a measurable functional
of $\left( \displaystyle\frac{1}{c}(\rho ^x_{c^2t})_{t\geq 0, x>0}\right)$ is the same
as the law of the rescaled functional of $\left( (\rho ^{x/c}_{t})_{t\geq 0, x>0}\right)$. 
In particular,
$$\mathcal L ((\int_0^t\frac{ds}{(\rho ^x_s)^2})_{t\geq 0, x>0}) = 
\mathcal L ((\int_0^{t/x^2}\frac{ds}{(\rho ^1_s)^2})_{t\geq 0, x>0}),$$
and we note using Cauchy sequence characterisation of P-a.s. convergence, that P-a.s.
convergence  for the rescaled processes is equivalent to
the same for the original processes.
\end{rem}

%%%%%%%%%%%%%%%%%%%%%%
\begin{lem}\label{spl}
If $\delta =2$ and $x\geq 0$, then as $t\rightarrow +\infty $,
$$\displaystyle\frac{4}{(\ln t)^2}\displaystyle\int _1^t
\displaystyle\frac{1}{(\rho ^x_s)^2} ds \stackrel{\mathcal L}{\rightarrow}
T_1(\beta)$$
where $T_1(\beta )$ is the first passage time of level 1 for a standard brownian motion. 

If  $\delta >2$ and $x\geq 0$, then as $t\rightarrow +\infty $,
$$\displaystyle\frac{1}{\ln t}\displaystyle\int _1^t
\displaystyle\frac{1}{(\rho ^x_s)^2} ds \stackrel{a.s.}{\rightarrow}
\displaystyle\frac{1}{(\delta - 2)}.$$

If  $\delta = 2$ and $x\geq 0$, then as $t\rightarrow +\infty $,
$$\displaystyle\frac{1}{\ln t}\displaystyle\int _1^t
\displaystyle\frac{1}{(\rho ^x_s)^2} ds \stackrel{a.s.}{\rightarrow}
+\infty.$$

\end{lem}
\vskip 0.5cm
\noindent{\it Proof\:} The first two asymptotics are well-known. For instance,
we can find the proof of the first one in \cite{MY} and, the second
and the third ones can be found in \cite{Ch}. 
$\Box $
\vskip 0.5cm

%%%%%%%%%%%%%%%%%%%%%%
\begin{lem}\label{spm}
If $\delta =2$ and $t>0$, then as $x\rightarrow 0+$,
$$\displaystyle\frac{1}{(\ln x)^2}\displaystyle\int _0^t
\displaystyle\frac{1}{(\rho ^x_s)^2} ds \stackrel{\mathcal L}{\rightarrow}
T_1(\beta)$$
where $T_1(\beta )$ is the first passage time of
 level 1 for a standard brownian motion. 

If  $\delta >2$ and $t>0$, then as $x\rightarrow 0+$,
$$\displaystyle\frac{1}{\ln x}\displaystyle\int _0^t
\displaystyle\frac{1}{(\rho ^x_s)^2} ds \stackrel{a.s.}{\rightarrow}
\displaystyle\frac{2}{(2 - \delta )}.$$

If  $\delta = 2$ and $t>0$, then as $x\rightarrow 0+$,
$$\displaystyle\frac{1}{\ln x}\displaystyle\int _0^t
\displaystyle\frac{1}{(\rho ^x_s)^2} ds \stackrel{a.s.}{\rightarrow}
-\infty.$$
The mentionned a.s. convergences are uniform in $t$ on compacts of $(0, +\infty )$.
\end{lem}
\vskip 0.5cm

\noindent{\it Proof\:} Take $x>0$ and consider the integral on $[0,t]$. 
Make a change of variables $s= x^2 u$ and 
use scaling property of Lemma \ref{sl}, Lemma \ref{spl} with $x=1$
and the Remark \ref{r1a}.
$\Box $
\vskip 0.5cm

%%%%%%%%%%%%%%%%%%%%%%%
\begin{lem}\label{inte}
 Let $\alpha \geq 0$ and $\beta \geq 0$. The integral
\begin{equation}\label{inte1}
\displaystyle\int_0^{+\infty }\displaystyle
\frac{(Y^1_s)^{\alpha }}{(\rho ^1_s)^{\beta }}ds
\end{equation}
is convergent P-a.s. iff
$\:\:\alpha(\delta -1) + (\beta-2)(\delta -2)>0.\:\:$
In particular,
it is convergent when $\alpha +\beta >2$.
The mentionned convergence is uniform in $t$ on compacts of 
$(0, +\infty )$.
\end{lem}
\vskip 0.5cm

\noindent{\it Proof\:}
We notice that the integrand in (\ref{inte1}) is positive
and, hence, the integral is
convergent P-a.s. to a finite limit or to $+\infty $.  
We remark that for $s\geq 0$ the first derivative of Bessel flow
$$Y^1_s =  \exp \{ - \displaystyle\frac{(\delta -1)}{2}
\displaystyle\int_0^s\displaystyle\frac{1}{(\rho ^1_u)^2}du\}.$$
By Ito formula we have:
$$\ln (\rho ^1_s) = 
\displaystyle\int_0^s\displaystyle\frac{1}{\rho ^1_u}d\beta _u +
\displaystyle\frac{(\delta -2)}{2}
\displaystyle\int_0^s\displaystyle\frac{1}{(\rho ^1_u)^2}du.$$
So, performing a time change  with $A_s =
\displaystyle\int_0^s\displaystyle\frac{1}{(\rho ^1_u)^2}du$
in considered integral we obtain that:
$$\displaystyle\int_0^{+\infty }\displaystyle
\frac{(Y^1_s)^{\alpha }}{(\rho ^1_s)^{\beta }}ds
\stackrel{\mathcal L}{=}
\displaystyle\int_0^{+\infty }\exp\left(a\tilde{\beta _u} - \frac{b\,u}{2}\right)du$$
where $a = 2-\beta $ and $b = \alpha(\delta -1) + (\beta-2)(\delta -2)$
and $\tilde{\beta }$ standard Brownian motion. It remains to note that 
the last integral is convergent iff $ b >0.$ $\Box $
\vskip 0.5cm

%%%%%%%%%%
To prove the existence of the derivatives of Bessel flow at $x=0$ we
need an explicite expression for the derivatives of Bessel flow at $x>0$.
For this we introduce $h^x = (h^x_t)_{t\geq 0}$ with
\begin{equation}\label{he}
h^x_t = - \displaystyle\frac{(\delta -1)}{2}
\displaystyle\int_0^t
\displaystyle\frac{1}{(\rho ^x_s)^2}ds. 
\end{equation}

Let now $\mathcal I_n$ be a set of multi-indices:
$$\mathcal I_n = 
\{ I = (i_1, i_2, \cdots i_n):
i_1\geq 0,i_2\geq 0,\cdots i_n\geq 0,\sum_{r=1}^n ri_r = n\}.$$
For $g\in C^n(\R)$ and $I \in \mathcal I _n$, $I = (i_1, i_2, \cdots i_n)$,
we introduce differential monomials
\begin{equation}\label{q}
Q_I(g) = 
(\frac{\partial g}{\partial x})^{i_1}
(\frac{\partial ^2 g}{\partial x^2})^{i_2}\cdots
(\frac{\partial ^n g}{\partial x^n})^{i_n}.
\end{equation}
as well as differential polynomials
\begin{equation}\label{pn}
P_n(g) = \sum _{I\in \mathcal I_n} 
c_IQ_I(g) 
\end{equation}
where $c_I$ are real constants.

Then using the existence of the first derivative
and recurrence arguments we obtain that for all $n\geq 1$
\begin{equation}
\label{dron}
\frac{\partial^n \rho ^x_{t\wedge \tau_{\epsilon} }}{\partial x^n} =
Y^x_{t\wedge \tau_{\epsilon}} P_{n-1}(h^x_{t\wedge \tau_{\epsilon}})
\end{equation}
with $P_0(\cdot ) = 1$ and $ \tau_{\epsilon}$ defined by (\ref{te}).
By the same reasoning we can prove that for all $1\leq k < n$
\begin{equation}\label{dhn}
\frac{\partial ^k h^x_{t\wedge \tau_{\epsilon}}}{\partial x^k} =
\sum _{I\in \mathcal I_k}
c_I\displaystyle\int_0^{t\wedge \tau_{\epsilon}}
\displaystyle\frac{Q_I(\rho ^x_s)}{(\rho ^x_s)^{2+j_k}}ds. 
\end{equation}
where $Q_I(\cdot)$ are the differential monomials of the type (\ref{q})
and for $I = (i_1, i_2, \cdots i_k)$, $j_k = \sum _{r=1}^k i_r$.
Taking the limit as $\epsilon \rightarrow 0 $ and using the fact
that $P(\tau_0 = \infty) =1$ we obtain the needed formulas. These formulas
coincide with the ones obtained by replacing of  $t\wedge \tau _{\epsilon}$ by $t$ in
(\ref{dron}) and(\ref{dhn}).

Now for $n\geq 1$ we introduce the integrals
\begin{equation}
\label{bn}
b_n=x^n\int _0^{\infty}\frac{(Y^x_s)^n}{(\rho ^x_s)^{n+2}}ds
\end{equation}
which are convergent according to Lemma \ref{inte}. Moreover, using
change of variables $s=s'x^2$ we establish that
$$b_n\stackrel{\mathcal L}{=}\int _0^{\infty}\frac{(Y^1_s)^n}{(\rho ^1_s)^{n+2}}ds$$
and that the law of $b_n$ does'not depend on $x$.
Let $B_0=1$ and
let us denote by $B_n$ the quantity:
\begin{equation}
\label{Bn}
B_n= \sum _{k=1}^n b_k^{1/k}
\end{equation}
We notice that the law of $B_n$ does'not depend on $x$ since it is so
for the $b_k$.
%%%%%%%%%%%%%%%%%%%%%
\begin{lem}\label{maj}
For each $n\geq 1$ there exists a real positive constant $c=c(n)$ such
that
$$|\,x^{n-1} \,\displaystyle\frac{\partial ^n \rho^x_t}{\partial x^n}\,|\leq
c\, Y^x_t\, B_{n-1}^{n-1},\,\,\
,|\,x^n \,\displaystyle\frac{\partial ^n h^x_t}{\partial x^n}\,|\leq c\,
B_n^n.$$
As a consequence for each $n\geq 1$ there exists a constant $c=c(n)$
such that 
$$|\,x^n\, P_n(h^x_t)\,| \leq c\, B_n^n. $$
\end{lem}

\noindent{\it Proof\:} The proof is going by induction using previous
formulas for the derivatives. For $n=1$ we have  
\begin{equation}\label{h11}
x\,\displaystyle\frac{\partial }{\partial x}(h^x_t) =
(\delta -1)\,x
\displaystyle\int_0^t\displaystyle\frac{Y^x_s}{(\rho ^x_s)^3}ds 
\end{equation}
and, hence, 
$$|x\,\displaystyle\frac{\partial }{\partial x}(h^x_t)|\leq (\delta -1)\,b_1.$$
So, we see that the claim is true with $c= \max (1, \delta -1 )$.
Suppose that for $1\leq m\leq n$
$$|\,x^{m-1} \,\displaystyle\frac{\partial ^m \rho^x_t}{\partial x^m}\,|\leq
c\, Y^x_t\, B_{m-1}^{m-1},\,\,\
,|\,x^m \,\displaystyle\frac{\partial ^m h^x_t}{\partial x^m}\,|\leq c\,
B_m^m.$$
We show that the needed relations hold for $m=n+1$. Below we will denote by
$c$ a generic constant.
From the formula (\ref{dron}) with replacing $n$ by $n+1$ and
(\ref{pn}) 
it follows that for the first estimation it is sufficient to majorate
 each $Q_I(h^x_t)$ with $I\in \mathcal I_n$ where
$$Q_I(h^x_t) = 
(\frac{\partial h^x_t}{\partial x})^{i_1}
(\frac{\partial ^2 h^x_t}{\partial x^2})^{i_2}\cdots
(\frac{\partial ^n h^x_t}{\partial x^n})^{i_n}$$
and $I=(i_1,i_2,\cdots , i_n)$ with $\sum _{r=1}^n r\,i_r = n$.
Since $B_n$ is increasing sequence we obtain from our suppositions
that
$$|\,x^n\,Q_I(h^x_t)\,| \leq c B_1^{i_1}B_2^{2i_2}\cdots B_n^{ni_n}\leq c
B_n^n$$
and it gives the first and third estimations of Lemma.

We remark that 
$$|\,\frac{\partial^{n+1} h^x_t}{\partial x^{n+1}}\,|\leq
c \sum _{I\in \mathcal I_{n+1}} \int _0^{\infty }\frac{|Q_I(\rho^x_s)|}{(\rho
^x_s)^{j_{n+1}+2}}ds$$
with $I=(i_1,i_2,\cdots ,i_{n+1})$,  $\sum _{r=1}^{n+1} r\,i_r = n+1$
and $j_{n+1}=\sum_{r=1}^{n+1} i_r$.

 Then we take  in account the formula
(\ref{dron}) to obtain
$$Q_I(\rho^x_s) =
(Y^x_s)^{j_{n+1}}(P_1(h^x_s))^{i_2}\cdots (P_n(h^x_s))^{i_{n+1}}.$$
Since $\sum _{r=2}^{n+1}( r-1)\,i_r = n+1 - j_{n+1}$, we have from previous estimations that
$$|\,x^{n+1-j_{n+1}}Q_I(\rho_s^x)\,|\leq c (Y^x_s)^{j_{n+1}} B_1^{i_2}B_2^{2i_3}\cdots B_n^{ni_{n+1}}.$$
Using this estimation and doing the change of variables in the integrals we
obtain that
$$|\,x^{n+1} \displaystyle\frac{\partial ^{n+1} h^x_t}{\partial
  x^{n+1}}\,|\leq c \sum _{I\in \mathcal I_{n+1}}b_{j_{n+1}} B_1^{i_2}B_2^{2i_3}\cdots
B_n^{ni_{n+1}}\leq c B_{n+1}^{n+1}$$
since $1\leq j_{n+1}\leq n+1$ and  $b_r\leq B_r^r$ for each $r$. So, we have the second estimation
and it proves Lemma.
$\Box $

%%%%%%%%%%%%%%%%%%%%%
\begin{lem}\label{cv}
Let $t>0$ and $x>0$. For all fixed $n\geq 1$ the sequences of random variables
$(x^n \displaystyle\frac{\partial ^n h^x_t}{\partial x^n})$ and
$(x^n P_n(h^x_t))$ are convergent P-a.s. as $x\rightarrow 0+$.
The mentionned convergences are uniform in $t$ on compacts of 
$(0,+\infty )$.
\end{lem}
\vskip 0.5cm

\noindent{\it Proof\:} We prove in recurrent way that for all $n\geq 1$
the sequence of random variables
 $(x^n \displaystyle\frac{\partial ^n h^x_t}{\partial x^n})$
is convergent  P-a.s. as $x\rightarrow 0+$ uniformly  in $t$ on compacts of 
$(0,+\infty )$. We remark that by formulas (\ref{q}), (\ref{pn}) 
this convergence gives immediately the same type of convergence 
for $(x^n P_n(h^x_t))$.

For $n=1$ we have
\begin{equation}\label{h1}
x\displaystyle\frac{\partial }{\partial x}(h^x_t) =
(\delta -1)x
\displaystyle\int_0^t\displaystyle\frac{Y^x_s}{(\rho ^x_s)^3}ds 
\stackrel{\mathcal L}{=}
(\delta -1)
\displaystyle\int_0^{t/x^2}\displaystyle\frac{Y^1_s}{(\rho ^1_s)^3}ds.
\end{equation}
By Lemma \ref{inte} with $\alpha = 1$ and $\beta = 3$ we prove that
the integral is convergent P-a.s., uniformly in $t$ on compacts of 
$(0,+\infty )$. 

Suppose that P-a.s. convergence, uniform  in $t$ on compacts of 
$(0,+\infty )$ is valid for 
$(x^k \displaystyle\frac{\partial ^k h^x_t}{\partial x^k})$
with $1\leq k\leq n$. Then by (\ref{pn}) we obtain that the
$(x^k P_k(h^x_t))$ are convergent P-a.s., uniformly  in $t$ on compacts of 
$(0,+\infty )$, as $x\rightarrow 0+$.
By formula (\ref{dhn}) we have:
$$\frac{\partial ^{n+1} h^x_t}{\partial x^{n+1}} =
\sum _{I\in \mathcal I_{n+1}}
c_I\displaystyle\int_0^t
\displaystyle\frac{Q_I(\rho ^x_s)}{(\rho ^x_s)^{2 + j_{n+1}}}ds.$$
where  $I=(i_1,i_2,\cdots ,i_{n+1})$,  $\sum _{r=1}^{n+1} r\,i_r = n+1$
and $j_{n+1}=\sum_{r=1}^{n+1} i_r$. We show that each term in the
previous sum multiplying by $x^{n+1}$ is convergent P-a.s. uniformly
in $t$ on compacts of $(0, \infty )$. For this we remark that the
term corresponding to $I = (i_1, i_2, \cdots i_{n+1})$ in the previous sum
times $x^{n+1}$ is equal to
$$x^{n+1}\,\displaystyle\int_0^t\displaystyle\frac{(Y^x_s)^{j_{n+1}}(P_1(h^x_s))^{i_2}
(P_2(h^x_s))^{i_3}\cdots (P_n(h^x_s))^{i_{n+1}}}{(\rho ^x_s)^{2 + j_{n+1}}}ds$$
and that  it is equal in law to
$$\displaystyle\int_0^{t/x^2}
\displaystyle\frac{(Y^1_s)^{j_{n+1}}(P_1(h^1_{s/x^2}))^{i_2}
(P_2(h^1_{s/x^2}))^{i_3}\cdots (P_n(h^1_{s/x^2}))^{i_{n+1}}}
{(\rho ^1_s)^{2+j_{n+1} }}ds$$
where the $P_k(h^1_{s/x^2})$ are equal in law to the $x^k P_k(h^x_s)$.

We notice that for $1\leq k\leq n$ the $(x^k P_k(h^x_t))$ 
are uniformly bounded by $B_k^k$ according to Lemma \ref{maj} and that
the law of $B_k^k$ does'not depend on $x$. Moreover,
since $j_{n+1} \geq 1$, the  integral

$$\displaystyle\int_0^{\infty }\displaystyle\frac{(Y^x_s)^{j_{n+1}}}{(\rho ^x_s)^{2 + j_{n+1}}}ds$$
is converging. So, changing space, we have P-a.s. convergence by
Lebesgue dominated convergence theorem. Then, the final result follows  by 
Lemma \ref{sl}. 

The same can be done simultaneously for the expression of 
$(x^{n+1}\displaystyle\frac{\partial ^{n+1} h^x_t}{\partial x^{n+1}})$ and this
proves P-a.s.
convergence of this variable, uniform in $t$ on compacts of
$(0, +\infty )$.
$\Box $
\vspace{0.5cm}

\begin{lem}\label{p}
 For $n\geq 1$ and $t>0$ let $U_n= \displaystyle\lim_{x\rightarrow 0+} x^n P_n(h^x_t)$. 
Then $U_n\neq 0$(P-a.s.).
\end{lem}

\noindent{\it Proof\:} Writing the expression for (n+2)-th derivative
of $\rho _t ^y$  from (\ref{dron}) and derivating the same expression for 
(n+1)-th derivative of $\rho _t^y$ we get that: for $y>0$ and $t>0$
$$P_{n+1}(h^y_t) =\frac{\partial}{\partial y}(h^y_t) P_n(h^y_t) + 
\frac{\partial}{\partial y}(P_n(h^y_t)).$$
We notice that  for $n\geq 1$
$$\lim _{y\rightarrow 0+ } y^{n+1}P_{n+1}(h^y_t) = U_{n+1},\:
\lim _{y\rightarrow 0+ } y^n P_n(h^y_t) = U_n,\:
\lim _{y\rightarrow 0+ } y\frac{\partial}{\partial y}(h^y_t) = U_1$$
and we prove that
\begin{equation}
\label{const}
U_{n+1} = (U_1  - n)U_n.
\end{equation}
For this we take $x>0$ and $\alpha\in ]0,1[$ and we integrate the previous equality on the
interval $[\alpha x, x]$:
$$\int _{\alpha x}^x P_{n+1}(h^y_t)dy = \int _{\alpha x}^x
\frac{\partial}{\partial y}(h^y_t) P_n(h^y_t)dy + 
P_n(h^x_t) - P_n(h^{\alpha x}_t).$$
Then we estimate  each integral, we multiply the result by $x^n$,
we take the limit as $x\rightarrow 0+$ and we obtain (\ref{const}).

Finally, we have
\begin{equation}
\label{cn}
U_{n+1} = U_1(U_1 -1)(U_1-2)\cdots (U_1-n)
\end{equation}
with
 $$U_1= (\delta -1)
\displaystyle\int_0^{\infty }\displaystyle\frac{Y^1_s}{(\rho ^1_s)^3}ds.$$ 

We show that the random variable $U_1$ has a density. For this
we  make a change of variable as in Lemma \ref{inte}
with $A_s =\displaystyle \int _0^s \frac{1}{(\rho ^1_u)^2}ds$
to prove that
$$U_1 \stackrel{\mathcal L}{=} 
(\delta -1)\int _0^{\infty} \exp \left(\beta_u - \frac{(2\delta -3)}{2}u\right) du$$
where $\beta $ is Brownian motion. Using Dufresne equality
(see for instance \cite{YY}, p. 95) we have
\begin{equation}
\label{c1}
U_1\stackrel{\mathcal L}{=}\displaystyle\frac{2(\delta -1)}{Z_{\nu}}
\end{equation}
where $Z_{\nu }$ follows gamma law $\Gamma (\nu ,1)$ of index $\nu = (2\delta -3)$.

Since $U_1$ has a density, $P(U_1\in
\N)=0$ and for each $n\geq 1$,
$U_n \neq 0$ (P-a.s.). $\Box $
\vspace{0.5cm}
          
 %%%%%%%%%%%%%%%%%%%%%%%%%%%%%%%%%%%%%%%%%%%%%%%%%%%%%%%%%%%%% 
\noindent {\it Proof of Theorems \ref{t11}\,\,} We have from (\ref{bes})
$$\rho ^0_t =  \beta _t + \displaystyle\frac{(\delta -1)}{2}
\displaystyle\int_0^t\displaystyle\frac{1}{\rho ^0_s}ds.$$
For $x>0$ and $t>0$ we put
$$Z^x_t = \displaystyle\frac{\rho ^x_t - \rho ^0_t}{x},$$
and we remark that $Z^x_t$ satisfies a linear stochastic equation
with the solution given by:
$$Z^x_t = \exp \{ - \displaystyle\frac{(\delta -1)}{2}
\displaystyle\int_0^t\displaystyle\frac{1}{\rho ^x_s\rho ^0_s}ds\}.$$
The fact that for all $s\geq 0$, $\rho ^x_s \downarrow \rho ^0_s$ (P-a.s.) 
as $x\downarrow 0$ and the property: P-a.s. for all $t>0$
\begin{equation}\label{int}
\displaystyle\int_0^t\displaystyle\frac{1}{(\rho ^0_s)^2}ds 
= +\infty ,
\end{equation}
together with Lebesgue monotone convergence theorem give that P-a.s. for all $t>0$
$$\displaystyle\frac{\partial \rho ^x_t}{\partial x}\mid _{x=0} =
\displaystyle\lim _{x\rightarrow 0+} Z^x_t = 0.$$
In the same manner we establish that P-a.s. for all $t_0>0$
$$\displaystyle\lim _{\stackrel{x\rightarrow 0+}{t\rightarrow t_0}}
\displaystyle\frac{\partial \rho ^x_t}{\partial x} =
\displaystyle\lim _{\stackrel{x\rightarrow 0+}{t\rightarrow t_0}}
\exp \{ - \displaystyle\frac{(\delta -1)}{2}
\displaystyle\int_0^t\displaystyle\frac{1}{(\rho ^x_s)^2}ds \} = 0
$$
 and this proves the continuity of the first derivative.

To prove the existence of the $n$-th derivative at $x=0$ equal to zero
we show that there exists $n(\delta )$ such that for 
$2\leq n < n(\delta )$ (P-a.s.)
\begin{equation}\label{exdn}
\displaystyle\lim _{x\rightarrow 0+}\displaystyle\frac{1}{x}
\displaystyle\frac{\partial ^{n-1} \rho ^x_t}{\partial x ^{n-1}} =0.
\end{equation}
To find $n(\delta )$ we write that 
$$\displaystyle \frac{\partial ^{n-1} \rho ^x_t}{\partial x ^{n-1}}
= Y^x_t P_{n-2}(h^x_t)$$
where $h^x$ is defined by (\ref{he}) and  $P_{n}$ is given by (\ref{pn}).
From Lemma \ref{cv} we have that uniformly in $t$ on compact sets of
$(0,+\infty )$
$$\lim _{x\rightarrow 0+} x^{n-2} P_{n-2}(h^x_t) = U_{n-2}$$
where $U_{n-2}$ is different from zero with probability 1.
This means that $(\ref{exdn})$ is equivalent to
\begin{equation}
\label{equi}
\lim_{x\rightarrow 0+}\exp\left\{-\left[ \displaystyle\frac{(\delta -1)}{2}
\displaystyle\int_0^t\displaystyle\frac{1}{(\rho ^x_s)^2}ds +
(n-1)\ln x \right]\right\}=0
\end{equation}
If $\delta =2$ then applying Lemma  \ref{spm} we see that the last
relation holds for all $n\geq 2$ and we can put $n(\delta) = +\infty $. If $\delta >2$ then it is easy to
see that for $n(\delta )= 1+\displaystyle\frac{\delta -1}{\delta -2}$
and $n<n(\delta )$ the relation (\ref{equi}) holds and for $n>n(\delta )$
it fails. If $\delta >2$ and $n=n(\delta )$ then by Ito formula
$$\ln (\rho ^x_t) = \ln x + 
\displaystyle\int_0^t\displaystyle\frac{1}{\rho ^x_s}d\beta _s +
\displaystyle\frac{(\delta -2)}{2}
\displaystyle\int_0^t\displaystyle\frac{1}{(\rho ^x_s)^2}ds.$$
We apply a central limit theorem (see \cite{LSh} , p.472 ) for the martingale 
for $M = (M_t)_{t\geq 0}$ with
$$M_t = \displaystyle\frac{1}{\sqrt{|\ln (x)|}}
\displaystyle\int_0^{t/x^2}\displaystyle\frac{1}{\rho ^1_s}d\beta _s$$
obtained from original one by time change, to prove via Skorohod
representation theorem that the quantity appearing as the power in
exponential in (\ref{equi}), namely
$$ \displaystyle\frac{(\delta -1)}{2}\displaystyle\int_0^t
\displaystyle\frac{1}{(\rho ^x_s)^2}ds +
(n(\delta )-1)\ln x $$
behaves as 
$c\, \xi |\ln x|^{1/2}$
as $x\rightarrow 0$ where $\xi $ is standard $\mathcal N(0,1)$ random
variable and $c$ is some positive constant.
 Hence, (\ref{equi}) fails on the set $\{\xi >0 \}$ of
probability 1/2, as well as (\ref{exdn}).
 
For the bicontinuity of the $n$-th derivative at $x=0$ for 
$2\leq n < n(\delta )$
we prove that P-a.s. for all $t_0 > 0$
\begin{equation}\label{cdn}
\displaystyle\lim _{\stackrel{x\rightarrow 0+}{t\rightarrow t_0}}
\displaystyle\frac{\partial ^n \rho ^x_t}{\partial x ^n} =0.
\end{equation}
The proof of (\ref{cdn}) is going in the same way as (\ref{exdn}) using 
the fact that  the convergences
in Lemmas \ref{cv} and \ref{spm} are uniform in $t$ on compact sets of
$(0, +\infty )$.

The asymptotic relations a), a') b), c) follows from (\ref{dron}) and  Lemmas \ref{spm} and
\ref{p}. To prove d) we remark that according to Lemma \ref{maj} we
have
$$\sup_{\epsilon\leq t\leq T} |x^{n-1} P_{n-1}(h^x_t)| \leq c B_{n-1}^{n-1}$$
with $B_n$ defined by (\ref{Bn}). It remains only to show that
$B_{n-1}^{\gamma (n-1)}$ is integrable. Since  for $a,b\in \R^+$ and $\gamma >0$,\,
$(a+b)^{\gamma}\leq c(a^{\gamma}+b^{\gamma})$ with some constant $c$,
$B_{n-1}^{(n-1)\gamma}$ is integrable if $b_k^{\gamma (n-1)/k}$ are integrables for
$1\leq k\leq n-1$. Making time change  with $A_t=k^2\displaystyle
\int_0^t\displaystyle\frac{1}{(\rho ^1_s)^2}ds$
and using Dufresne identity
(see \cite{BS}, p.78) we obtain that
$$b_k\stackrel{\mathcal L}{=}\frac{1}{k^2}\int _0^{\infty }\exp (\beta_u
-\frac{(2\delta -3)u}{2k})\,du\stackrel{\mathcal L}{=}\frac{2}{kZ_{\nu (k)}}
$$
where $(\beta_u)_{u\geq0}$ is standard Brownian motion and $Z_{\nu (k)}$ 
is the variable following gamma law $\Gamma (\nu(k), 1)$
of index $\nu (k)=\frac{2\delta -3}{k}$.
 So, we have needed integrability if for all
$1\leq k\leq n-1$, the variables $(2/Z_{\nu (k)})^{\frac{\gamma(n-1)}{k}}$
are integrables. As well-known this is true, if
$\frac{\gamma(n-1)}{k}< \nu (k)$ and
the last condition is satisfied for $\gamma < \frac{2\delta -3}{n-1}$.
$\Box $
\vskip 0.5cm

%%%%%%%%%%%%%%%%%%%%%%%%%%%%%%%%%%%%%%%%%%%%%%%%%%%%%%%%%%%%%%%%%%%%%%%%%%%%%%%%%%%%%%%%%%%%
            \section{Regularity of Bessel flow for $1 < \delta < 2$}
%%%%%%%%%%%%%%%%%%%%%%%%%%%%%%%%%%%%%%%%%%%%%%%%%%%%%%%%%%%%%%%%%%%%%%%%%%%%%%%%%%%%%%%%%%%%
 \vskip 0.5cm

We start with some lemmas needed to prove the theorem \ref{t2}.
Let $\tau_0(x)$ be defined by (\ref{tauo}).

%%%%%%%%%%%%%%%%%%%%%%%  
\begin{lem}\label{ltau}
Let $x>0$ be fixed. Then $\tau _0(x-)=\tau _0(x)=\tau _0(x+)$ (P-a.s.). 
Moreover, there exists a cadlag version of $(\tau _0(x))_{x>0}$.
\end{lem}
\vskip 0.5cm

\noindent {\it Proof}
By comparison theorem we have that for $x\leq y$ (P-a.s.)
$$\tau _0 (x) \leq \tau _0 (y).$$ Then, there exist the limits (P-a.s.):
$$\lim_{y\rightarrow x+}\tau _0(x) = \tau _0(x+), \:
\lim_{y\rightarrow x-}\tau _0(x) = \tau _0(x-).$$
We take $y<x<z$ then by comparison theorem again for $\gamma >0$
$$E(\tau _0(y)^{\gamma})\leq E(\tau _0(x)^{\gamma})\leq E(\tau _0(z)^{\gamma})$$
Since
\begin{equation}\label{id}
\tau _0(x)\stackrel{\mathcal L}{=} \displaystyle\frac{x^2}{2\gamma _{\nu}},
\end{equation}
where $\gamma _{\nu}$ is random variable of gamma law $\Gamma (\nu
,1)$ with index 
$\nu = 1 - \frac{\delta}{2}$, we have for small $\gamma >0$ that
$$E(\tau _0(x))^{\gamma} = c_{\gamma } x^{2\gamma}$$
with some positive constant  $c_{\gamma }$.
The mentionned informations implies that for small $\gamma >0$
$$E(\tau _0(x-)^{\gamma})= E(\tau _0(x)^{\gamma})= E(\tau
_0(x+)^{\gamma})$$
and, hence, P-a.s. $\tau _0(x-)=\tau _0(x)=\tau _0(x+)$.

To construct a cadlag version of $(\tau _0(x))_{x>0}$, we take for $x\in \Q^+$ the value of $\tau _0(x)$
and for $x\in \R^+\setminus \Q^+$ we put:
$$\tau _0(x-) = \lim _{\stackrel{y\rightarrow x-}{y\in \Q^+}}\tau _0(y),\,\,
\tau _0(x) = \lim _{\stackrel{y\rightarrow x+}{y\in \Q^+}}\tau _0(y).$$
This construction gives a cadlag version of $(\tau _0(x))_{x>0}$ which preserves the
finite-dimensional distributions of $\tau _0(x)$. 
$\Box $
\vskip 0.5cm

%%%%%%%%%%%%%%%%%%%%%%%%
\begin{lem}\label{first}
We consider bicontinuous (P-a.s.) modifications of $(\rho^x_t)_{t\geq 0, x>0}$.
Then the flow of processes $(\rho^x_{t\wedge\tau _0(x)})_{t\geq 0}$ is bicontinuous in probability
and has the  first derivative in probability sense
with respect to $x$. This derivative is bicontinuous in probability on the set 
$]0, +\infty [\times ]0, +\infty ]$.
\end{lem}
\vskip 0.5cm

\noindent {\it Proof }Let $(x,t)$ be fixed with $x>0, t\geq 0$. To prove bicontinuity  in probability
of $(\rho^x_{t\wedge\tau _0(x)})_{x>0, t\geq 0}$, we write:
$$|\rho^y_{s\wedge\tau _0(y)} - \rho^x_{t\wedge\tau _0(x)}|\leq
|\rho^y_{s\wedge\tau _0(y)} - \rho^x_{s\wedge\tau _0(y)}| +
|\rho^x_{s\wedge\tau _0(y)} - \rho^x_{t\wedge\tau _0(x)}|.$$
We introduce the set
$A_{\delta } = \{|\tau (x) - \tau (y)| \leq \delta \}$
with $0<\delta <1$.
Then, on the set $A_{\delta}$ for $|s-t| < \gamma < \delta $ we have:
$$|\rho^y_{s\wedge\tau _0(y)} - \rho^x_{t\wedge\tau _0(x)}|\leq
\sup _{0\leq u\leq t+\delta}|\rho^y_u - \rho^x_u| + 
\sup _{|u-v|\leq \delta ;\, u,v \leq t+\delta}|\rho^x_u - \rho^x_v|.$$
We write for $\epsilon >0$ that
$$P(|\rho^y_{s\wedge\tau _0(y)} - \rho^x_{t\wedge\tau _0(x)}|\geq \epsilon ) \leq
P(\{|\rho^y_{s\wedge\tau _0(y)} - \rho^x_{t\wedge\tau _0(x)}|\geq \epsilon \}\cap A_{\delta }) +
P(A^c_{\delta }).$$
From the previous estimations we obtain that
$$P(|\rho^y_{s\wedge\tau _0(y)} - \rho^x_{t\wedge\tau _0(x)}|\geq \epsilon ) \leq
P(\sup _{0\leq u\leq t+1}|\rho^y_u - \rho^x_u|\geq \epsilon /2)  + 
P(\sup _{|u-v|\leq \delta ;\, u,v \leq t+1}|\rho^x_u - \rho^x_v|\geq \epsilon /2) +
P(A^c_{\delta }).$$
Using the facts that  $(\rho^x_u)_{ u\geq 0, x>0}$ is bicontinuous
P-a.s. and that $\tau (x)$
is continuous in probability, we obtain
taking $\lim _{\delta \rightarrow 0}\lim
_{y \rightarrow x}$, the claimed bicontinuity.

We show that for $t\in ]0, +\infty ]$ the first derivative of $(\rho^x_{t\wedge \tau _0(x)})_{t\geq 0, x>0}$ 
with respect to $x$ is given by:
\begin{equation}\label{6a}
Y^x_{t\wedge\tau _0(x)} = \exp \{ -\frac{(\delta -1)}{2}\int _0^{t\wedge\tau _0(x)}
\frac{ds}{(\rho^x_s)^2}\}
\end{equation}
For this we take $\epsilon >0$ and we do a localisation with
$$\tau _{\epsilon}(x) = \inf \{ s \geq 0 : \rho ^x_s \leq \epsilon\}.$$
Then we write (\ref{incr3}). We notice that 
$\tau _0(x) = \lim _{\epsilon \rightarrow 0}\tau _{\epsilon}(x)$
and, then,
\begin{equation}\label{9} 
Z^{x,y}_{t\wedge\tau _0(x)} = 
\exp \{ -\frac{(\delta -1)}{2}\int _0^{t\wedge\tau _0(x)}
\frac{ds}{\rho^x_s\rho^y_s}\}.
\end{equation}
where $Z^{x,y}_t$ is defined by (\ref{incr1}).
Via comparison theorem it can be shown that
$$\lim _{y\rightarrow x} \int _0^{t\wedge\tau _0(x)}
\frac{ds}{\rho^x_s\rho^y_s} = \int _0^{t\wedge\tau _0(x)}
\frac{ds}{(\rho^x_s)^2}$$

and taking $\lim_{y\rightarrow x}$ in (\ref{9}) we have (\ref{6a}). To obtain the result for $t=\infty$
it is sufficient to take $\displaystyle\lim _{t\rightarrow \infty }$ in (\ref{incr3}) and continue 
in above way.

By time reversal we have:
\begin{equation}\label{9a}
\mathcal L\left((\rho^{x, \delta}_{\tau _0(x)-t})_{0\leq t\leq \tau _0(x)}\right) = 
\mathcal L \left((\rho^{0, 4-\delta}_t)_{0\leq t\leq L(x)}\right)
\end{equation}
where $\rho^{x, \delta}$ is $\mbox{BES}^x(\delta )$ process and
$L(x) = \sup\{ t\geq 0: \rho^{0, 4- \delta}_t = x\}$. Since
$1<\delta<2$, we have that
$4-\delta >2$ and using asymptotics of Lemma \ref{spm} we obtain: (P-a.s.)
\begin{equation}\label{10}
\int _0^{\tau _0(x)}\frac{ds}{(\rho^x_s)^2} = +\infty
\end{equation}
and it gives together with (\ref{6a}) the expression
$$Y^x_{t\wedge\tau _0(x)} = \left\{\begin{array}{ccc}
\exp \{ -\frac{(\delta -1)}{2}\int _0^{t\wedge\tau _0(x)}\frac{ds}{(\rho^x_s)^2}\}
&\mbox{if} & t< \tau _0(x),\\\\
0&\mbox{if}& t\geq \tau _0(x).
\end{array}
\right.$$

Now we prove a bicontinuity of the first derivative at each point
$(x,t)$  with $x>0, t>0$.
 We consider tree sets $D_1, D_2,D_3$:
\begin{equation}\label{6c}
D_1 = \{\omega : \tau _0(x) >t \},\:\:
D_2 = \{\omega : \tau _0(x) <t \},\:\:
D_3 = \{\omega : \tau _0(x) =t \}
\end{equation}
and we prove a bicontinuity on each of them.
For $D_1$ we write that
$D_1 = \bigcup _{\epsilon >0} D_1^{\epsilon}$ where
$D_1^{\epsilon} = \{\tau_{\epsilon}(x) > t\}$. On each set $D_1^{\epsilon}$
bicontinuity of $Y^x_{t\wedge\tau _0(x)}$ follows from bicontinuity (P-a.s.) of $(\rho^x_t)_{t>0, x>0}$.
Hence, taking a countable set of $\epsilon $ we obtain the same result on $D_1$.
On the set  $D_2$ we have $Y^x_{t\wedge\tau _0(x)}= Y^x_{\tau _0(x)} = 0$, 
and, hence it is continuous.

Take now the set $D_3$, then $t=\tau _0(x)$. Let $(s,y)$ be in the neighbourhood of $(\tau _0(x), x)$.
We show that $Y^y_{s\wedge\tau _0(y)}$ is in the neighbourhood of $Y^x_{\tau _0(x)} =0$.
In fact, on the set $\{s\geq \tau _0(y)\}$, $Y^y_{s\wedge\tau _0(y)} = Y^y_{\tau _0(y)} = 0$.
On the set $\{s<\tau _0(y)\}$ we remark that for all $x>0$ (P-a.s.)
\begin{equation}\label{11}
\lim _{\stackrel{y\rightarrow x}{s\rightarrow \tau _0(x)-}}\int _0^s\frac{du}{(\rho^y_u)^2} = +\infty
\end{equation}
In fact, by comparaison theorem for small $\gamma >0$ and $\epsilon
>0$ we have
\begin{equation}\label{12}
\lim _{\stackrel{y\rightarrow x}{ s\rightarrow \tau _0(x)-}}\int _0^s\frac{du}{(\rho^y_u)^2} 
\geq \lim _{y\rightarrow x}\int _0^{(\tau _0(x)-\gamma)\wedge \tau _{\epsilon}(x)}
\frac{du}{(\rho^y_u)^2} = 
 \int _0^{(\tau _0(x)-\gamma)\wedge \tau _{\epsilon}(x)}
\frac{du}{(\rho^x_u)^2}
\end{equation}
Taking $\epsilon \rightarrow 0$ and then $\gamma \rightarrow 0$ we have from 
(\ref{10}) and (\ref{12}) the relation (\ref{11}).
$\Box $
\vskip 0.5cm
To investigate the existence of the derivatives of higher order, we 
prove the following two lemmas.

%%%%%%%%%%%%%%%%%%%%%%%%%
\begin{lem}\label{regt} 
Let $x>0$ be fixed. Then 
$$P\lim _{y\rightarrow x} \frac{\tau _0(y) - \tau _0(x)}{(y-x)^2} = 0.$$
\end{lem}
\vskip 0.5cm

\noindent {\it Proof} Let $y>x>0$.
Since $(\rho^z_{\tau _0(x)}, z>0)$ is $\mathcal F _{\tau _0(x)}$ - Markov, we have
\begin{equation}\label{12d}
\tau _0(y) - \tau _0(x) = \inf \{ u>0: \rho _u^{\rho^y_{\tau _0(x)}} \} = 
\tau _0(\rho^y_{\tau _0(x)})
\end{equation}
Using (\ref{id}) we obtain that
\begin{equation}\label{12e}
\mathcal L(\tau _0(\rho^y_{\tau _0(x)})) = 
\mathcal L((\rho^y_{\tau _0(x)})^2 \frac{1}{2\gamma _{\nu}}).
\end{equation}
where $\gamma _{\nu}$ is gamma random variable with index
$\nu = 1- \delta /2$, independent from $\rho^y_{\tau _0(x)}$.
But 
\begin{equation}\label{12f}
\frac{(\rho^y_{\tau _0(x)})^2}{(y-x)^2} = 
\frac{(\rho^y_{\tau _0(x)}-\rho^x_{\tau _0(x)})^2}{(y-x)^2}
\end{equation}
since $\rho^x_{\tau _0(x)} = 0$. 
It was shown in Lemma \ref{first} that (P-a.s.)
$$\lim _{y\rightarrow x} \frac{\rho^y_{\tau _0(x)} - \rho^x_{\tau _0(x)}}{y-x} = 0.$$
Then, (\ref{12d}), (\ref{12e}) and (\ref{12f}) implies 
$$P\lim _{y\rightarrow x+} \frac{\tau _0(y) - \tau _0(x)}{(y-x)^2} = 0.$$
The same consideration with $x>y>0$ gives again (\ref{12f}) with the exchanging $x$ and $y$.
But
$$\frac{\rho^x_{\tau _0(y)}-\rho^y_{\tau _0(y)}}{x-y}= 
Z^{x,y}_{\tau _0(y)} = 
\exp \{ -\frac{(\delta -1)}{2}\int _0^{\tau _0(y)}
\frac{ds}{\rho^x_s\rho^y_s}\},$$
and the relation $(P-a.s.)$
$$\lim _{y\rightarrow x-}\int _0^{\tau _0(y)} \frac{ds}{\rho^x_s\rho^y_s}\rightarrow +\infty$$
implies
$$P\lim _{y\rightarrow x-} \frac{\tau _0(y) - \tau _0(x)}{(y-x)^2} = 0$$
and it proves the result.
$\Box $
\vskip 0.5cm

%%%%%%%%%%%%%%%%%%%%%%
\begin{lem}\label{dho}
The flow of the  processes $(Y^x_{t\wedge \tau _0(x)})_{t>0}$ defined
by (\ref{6a})  has bicontinuous derivatives
in probability sense  at $x>0$ only up to the order $n<n(\delta )$ where 
$n(\delta ) = \displaystyle\frac{\delta -1}{2 - \delta }$.
\end{lem}
\vskip 0.5cm

\noindent {\it Proof} First of all we remark that 
$D_1 = \bigcup _{\epsilon >0} D_1^{\epsilon}$
and on the sets $D_1^{\epsilon}$ the flow of 
$(Y^x_{t\wedge \tau _0(x)})_{t>0}$ with $x>0$ has the bicontinuous derivatives of all
orders. This follows from the fact that on $D_1^{\epsilon}$ this process
coincide with $(Y^x_{t\wedge \tau _{\epsilon}(x)})_{t>0,x>0}$ and we can use
the previous results for classical case. On the set $D_2$ the result is also
trivially true.

 Let $(x,t)\in D_3$ be fixed with $x>0$ and 
$t = \tau _0(x)$. Since $Y^y_{\tau_0(y)}=0$ we evidently have 
that 
$$P\lim_{y\rightarrow x-}\displaystyle \frac{\partial^n Y^y_{t\wedge\tau_0(y)}}
{\partial y^n} =0.$$
If we show that there exists $n(\delta )>0$ such that for $n<n(\delta )$
\begin{equation}
\label{partial}
P\lim_{y\rightarrow x+}\displaystyle \frac{\partial^n Y^y_{t\wedge\tau_0(y)}}
{\partial y^n} =0,
\end{equation}
then the  mentionned relations and continuity of the derivatives on $D_1$
and $D_2$ will imply that the flow of the processes $(Y^x_{t\wedge \tau _0(x)})_{t>0}$
has continuous derivatives  of the order $n<n(\delta )$.
We recall that for $y>x$ and $t=\tau_0(x)$
\begin{equation}
\label{r1}
\frac{\partial^n Y^y_{t\wedge \tau_0(y)}}{\partial
y^n}=Y^y_{\tau_0(x)}P_n(h^y_{\tau_0(x)})
\end{equation}
where $h^y$ and $P_n$ are defined by (\ref{he}),(\ref{pn}).

Let $u=y/x.$
Performing time change $s= s'x^2$ we obtain 
\begin{equation}
\label{r2}
Y^y_{\tau_0(x)}P_n(h^y_{\tau_0(x)})\stackrel{\mathcal L}{=}
Y^u_{\tau_0(1)}P_n(h^u_{\tau_0(1)})/x^n.
\end{equation}
First of all we investigate the behaviour of
$$Y^u_{\tau _0(1)}= \exp\left\{- \displaystyle\frac{(\delta -1)}{2}
\displaystyle\int_0^{\tau_0(1)}\displaystyle\frac{1}{(\rho ^u_s)^2}ds
\right\}.$$
Using time reversal we obtain that
\begin{equation}
\label{rev}
\mathcal L\left( \int_0^{\tau _0(1)}\frac{1}{(\rho^u_s)^2}ds
\,|\,\rho _{\tau _0(1)}^u = v\right)= 
\mathcal L\left(\int_{0}^{L_u(v)}\frac{1}{(\rho^{v,4-\delta }_s)^2}ds
\,\right)
\end{equation}
where $L_u(v) = \sup \{ s\geq 0: \rho^{v,4-\delta }_s = u\}$.
By time change we obtain that
\begin{equation}
\label{id1}  
\int_{0}^{L_u(v)}\frac{1}{(\rho^{v,4-\delta }_s)^2}ds
\stackrel{\mathcal L}{=}
 \int_{0}^{L_{u/v}(1)}\frac{1}{(\rho^{1,4-\delta }_s)^2}ds 
\end{equation}
where $L_{u/v}(1)$ is defined as previously with replacing of $u$ by
$u/v$ and $v$ by $1$.
Since as $a\rightarrow +\infty $
$$\frac{L_a(1)}{a^2}\stackrel{\mathcal L}{\rightarrow }\tau _{0}(1)$$
where $\tau_0(1)$ is the corresponding time of attending of zero, we
have that  
$$P\,\lim_{a\rightarrow +\infty}\frac{\ln L_a(1)}{2\ln a}= 1$$
Then using Lemma \ref{spl} we obtain that 
\begin{equation}
\label{asyv}
P\,\lim _{v\rightarrow 0+}  \frac{1}{|\ln (v)|}
\int_{0}^{L_{u/v}(1)}\frac{1}{(\rho _s^{1, 4-\delta })^2}ds = 
\frac{2}{2-\delta }.
\end{equation}
Since $\displaystyle\lim _{u\rightarrow 1+}\rho _{\tau _0(1)}^u =0$, we obtain
using standard arguments from
(\ref{asyv}), (\ref{rev}) and (\ref{id1}) that
\begin{equation}
\label{asyr}
P\,\lim _{u\rightarrow 1+}  \frac{1}{|\ln (\rho _{\tau _0(1)}^u)|}
\int_{0}^{\tau _0(1)}\frac{1}{(\rho^u_s)^2}ds = 
\frac{2}{2-\delta }.
\end{equation}

 If we appy a time reversal to $P_n(h^u_{\tau_0(1)})$ then we
 obtain
$$\mathcal L \left( P_n(h^u_{\tau_0(1)})|\rho^u_{\tau_0(1)}=v\right) =
\mathcal L\left(P_n(\,h^{v,4-\delta}_{L_u(v)}\,)\right).$$

From Lemma \ref{cv} and \ref{p} we obtain that 

$$P \lim _{v\rightarrow 0+}v^n\,P_n(\,h^{v,4-\delta}_{L_u(v)}\,)= U_n$$
where  $U_n$ is defined by
formula (\ref{cn}) with $U_1$ given by (\ref{c1}) and $\nu =
5-2\delta $.
So, by standard arguments we deduce that 
\begin{equation}
\label{47a}
P \lim _{u\rightarrow 1+}(\rho ^u _{\tau _0(1)})^n\, P_n(h^u_{\tau _0(1)}) = U_n
\end{equation}
Since $U_n\neq 0$ with probability 1,
 we obtain finally that
the relation (\ref{partial}) is equivalent to
\begin{equation}
\label{relat}
\lim_{u\rightarrow 1+}\exp\left\{-\left[ \displaystyle\frac{(\delta -1)}{2}
\displaystyle\int_0^{\tau_0(1)}\displaystyle\frac{1}{(\rho ^u_s)^2}ds
 + n \ln (\rho ^u_{\tau _0(1)}) \right]\right\}=0.
\end{equation}

Let $n(\delta ) =(\delta-1)/(2 -\delta )$. It is easy to see from
(\ref{asyr}) and (\ref{47a}) that if $n<n(\delta )$ then (\ref{relat}) holds and that
if $n>n(\delta )$ then (\ref{relat}) fails. For $n=n(\delta )$ it can
be shown using a central limit theorem for martingales and Skorohod
representation theorem like in the proof of theorems
\ref{t11} that the limit in (\ref{relat}) exists only on the
set of probability 1/2, and, that (\ref{relat}) fails, too.

$\Box $
\vskip 0.5cm

\noindent{\it Proof of theorem \ref{t2}}
Let $\rho^0$ be continuous version of $\mbox{BES}^0(\delta )$ process starting
from zero and $(\tau _0(x))_{x>0}$ be a cadlag version of the
corresponding process. For all $x>0$ and $t\geq 0$ we put:
\begin{equation}\label{13}
\tilde{\rho}^x_t = \rho ^x_{t\wedge \tau _0(x)}I_{[\![0, \tau _0(x)[\![} +
(\rho ^0_{t} - \rho ^0_{\tau _0(x)})I_{[\![\tau _0(x), +\infty [\![}.
\end{equation}
Using strong Markov property  we prove that the both processes $\rho $
and $\tilde{\rho }$ have the same finite-dimensional distributions.

We show that the trajectories of (\ref{13}) are in $D(\R^{+,*}, C(\R^{+,*}))$.
For this we remark that for each $x>0$, the process
$(\tilde{\rho} ^{x}_t)_{t\geq 0}$  is continuous in $t$. Moreover,
for $y>x$
$$|\tilde{\rho} ^y_t - \tilde{\rho} ^x_t| \leq
2\displaystyle\sup _{u,v \in U(x,y)}|\rho ^y_u - \rho ^x_v|+ 
\displaystyle\sup _{u,v\in U(x,y)}|\rho ^0_u - \rho
^0_v|$$
where $U(x,y)=\{(u,v):|u-v|\leq (\tau _0(y)- \tau _0(x))\,;\,
 u,v\leq t+(\tau _0(y)- \tau _0(x))\}$.
Since for all $x>0$ ,
$$\lim _{y\rightarrow x+}\tau _0(y) = \tau _0(x+) = \tau _0(x),$$
we see, that uniformly on compact sets of $t$, the right-hand side of the last
inequality is tending to zero as $y\rightarrow x+$.
Taking $y<x$ and using the fact that  for all $x>0$
$$\lim _{y\rightarrow x-}\tau _0(y) = \tau _0(x-),$$
we obtain the existence of left-hand limits uniformly on compact set
of $t$.

We show that  the first derivative of $\tilde{\rho}$ coinside with the one of
$(\rho ^x_{t\wedge\tau _0(x)})_{t\geq 0, x>0}$.
In fact, consider tree sets $D_1, D_2, D_3$ defined by (\ref{6c}). On the set
$D_1$ the process $(\tilde{\rho} ^{x}_t)_{t\geq 0}$ coincide with
$(\rho ^x_{t\wedge\tau _0(x)})_{t\geq 0, x>0},$ and the existence of
the first derivative was already
discussed in Lemmas \ref{first}, \ref{dho}. On $D_2$ the same process 
coincide with $(\rho ^0_{t})_{t\geq \tau _0(x), x>0}$ according to
comparison theorem
 and the first derivative 
 is equal to zero.
0n the set $D_3$ we have $t=\tau _0(x)$ and
$$\frac{\tilde{\rho} ^y_t - \tilde{\rho} ^x_t}{y-x} = \left\{
\begin{array}{ccc}
\displaystyle\frac{\rho ^y_{\tau _0(x)} - \rho ^x_{\tau _0(x)}}{y-x}&\mbox{if}&y>x,\\\\
0&\mbox{if}&y<x.
\end{array}\right.$$

We obtain as in Lemma \ref{regt} that
$$\lim _{y\rightarrow x+}\frac{\rho ^y_{\tau _0(x)} - \rho ^x_{\tau
    _0(x)}}{y-x}=0.$$
Hence, the first derivative of $\tilde{\rho}$ coinside with the one of
$(\rho ^x_{t\wedge\tau _0(x)})_{t\geq 0, x>0}$
and  we obtain the  claims
from Lemma \ref{dho}, and from the relations
(\ref{r1}),(\ref{r2}),(\ref{asyr}) and (\ref{47a}).

 Now we consider the case $x=0$ and $t>0$. 
We put $(\tilde{\rho}  ^0_t)_{t\geq 0} = (\rho ^0_t)_{t\geq 0}$ and we
remark that for $x>0$
\begin{equation}
\label{eq1}
 \tilde{\rho} ^x_t -  \tilde{\rho} ^0_t =  \rho ^x_{t\wedge
   \tau _0(x)} - \rho ^0_{t\wedge \tau _0(x)}
\end{equation}

Making time change and using scaling Lemma we have:
$$\frac{1}{x}( \rho ^x_{t\wedge
   \tau _0(x)} - \rho ^0_{t\wedge \tau _0(x)})
 \stackrel{\mathcal L}{=} \rho ^1_{(t/x^2)\wedge
   \tau _0(1)} - \rho ^0_{(t/x^2)\wedge \tau _0(1)}$$

Since $P(\tau _0(1)<\infty)=1$, we obtain that for $t>0$
$$P\lim _{x\rightarrow 0+}\frac{\tilde{\rho} ^x_t -  \tilde{\rho} ^0_t
}{x}=0.$$
Then we verify easily a bicontinuity of the first derivative at  $x=0$ and $t>0$.
For $x>0$ we have
$$Y^x_t = Y^x_{t\wedge \tau _0(x)} \stackrel{\mathcal
  L}{=}Y^1_{t/x^2\wedge \tau _0(1)}.$$
Since $P(\tau _0(1)<\infty)=1$ and (\ref{10}) we obtain that for each
$\epsilon >0$
$$P(Y^1_{t/x^2\wedge \tau _0(1)}>\epsilon) = P(\tau _0(1)>
  t/x^ 2) \rightarrow 0$$
as $x\rightarrow 0+$.

From the formula (\ref{dron}) we have for $x>0$ and 
$1\leq n < n(\delta )$ that
$$\frac{\partial^n \tilde{\rho }^x_t}{\partial x^n} =
\frac{\partial^{n-1} Y ^x_{t\wedge \tau_0(x) }}{\partial x^{n-1}} =
Y^x_{t\wedge \tau_0(x)} P_{n-1}(h^x_{t\wedge \tau_0(x)})$$
and, hence,
$$\frac{\partial^n \tilde{\rho }^x_t}{\partial x^n} \stackrel{\mathcal L}{=}
Y^1_{(t/x^2)\wedge \tau_0(1)} P_{n-1}(h^1_{(t/x^2)\wedge
  \tau_0(1)})/x^{n-1}.$$
We can see from the Lemma \ref{inte} and the proof of Lemma \ref{cv}
that for $\delta > 3/2$, $P_{n-1}(h^1_{\tau_0(1)})$ is a finite
random variable. So, for each $\epsilon >0$
$$P(\frac{1}{x}\frac{\partial^n \tilde{\rho }^x_t}{\partial
  x^n}>\epsilon)= P(\tau _0(1)>t/x^2)\rightarrow 0$$
as $x\rightarrow 0+$. It means that for $3/2<\delta <2$ there exist
the bicontinuous derivatives of  order $1\leq n < n(\delta)$ with
respect to $x$ in probability sense at $x\geq 0$ and $t>0$.
$\Box $

%%%%%%%%%%%%%%%%%%%%%%%%%%%%%%%%%%%%%%%%%%%%%%%%%%%%%%%%%%%%%%%% 
 \section{Acknowledgement}
 %%%%%%%%%%%%%%%%%%%%%%%%%%%%%%%%%%%%%%%%%%%%%%%%%%%%%%%%%%%%%%%%%
 The author is very grateful to Marc Yor
 for his helpful suggestions  and discussions on this topic.

%%%%%%%%%%%%%%%%%%%%%%%%%%%%%%%%%%%%%%%%%%%%%%%%%%%%%%%%%%%%%%%%%%%%%%%%%%%%%%%%
%                               bibliographie                                  
%%%%%%%%%%%%%%%%%%%%%%%%%%%%%%%%%%%%%%%%%%%%%%%%%%%%%%%%%%%%%%%%%%%%%%%%%%%%%%%%

%%%%%%%%%%%%%%%%%%%%%%%%%%%%%%%%%%%%%%%%%%%%%%%%%%%%%%%%%%%%%%%%%%%%%%%%%%%%%%%%%  
\end{document}